\magnification\magstephalf
\baselineskip14pt
\parskip2pt

\def\bn{\bigskip\noindent}
\def\ra{\rightarrow}
\def\bib{\par\noindent\hangindent 20pt}

\chardef\other=12
\outer\def\begintt{$$\let\par=\endgraf \ttverbatim \parskip=0pt
  \rightskip-2\parindent \ttfinish}
{\catcode`\|=0 |catcode`|\=\other 
  |obeylines 
  |gdef|ttfinish#1^^M#2\endtt{#1|vbox{#2}|endgroup$$}}
\def\ttverbatim{\begingroup \catcode`\@=\other \catcode`\"=\other
  \catcode`\\=\other
  \catcode`\{=\other
  \catcode`\}=\other
  \catcode`\$=\other
  \catcode`\&=\other
  \catcode`\#=\other
  \catcode`\%=\other
  \catcode`\~=\other
  \catcode`\_=\other
  \catcode`\^=\other
  \obeyspaces \obeylines \tt}
{\obeyspaces\gdef {\leavevmode\space}}

\centerline{\bf Convolution Polynomials}
\bigskip
\centerline{Donald E. Knuth}
\centerline{Computer Science Department}
\centerline{Stanford, California 94305--2140}

\bigskip
{\narrower\smallskip\noindent
{\bf Abstract.} The polynomials that arise as coefficients when a
power series is raised to the power~$x$ include many important special
cases, which have surprising properties that are not widely known.
This paper explains how to recognize and use such properties, and it
closes with a general result about approximating such polynomials
asymptotically.
\smallskip}

\bn
A family of polynomials $F_0(x),F_1(x),F_2(x),\ldots$ forms a {\it
convolution family\/} if $F_n(x)$ has degree $\leq n$ and if the
convolution condition
$$F_n(x+y)=F_n(x)F_0(y)+F_{n-1}(x)F_1(y)+\cdots
+F_1(x)F_{n-1}(y)+F_0(x)F_n(y)$$
holds for all $x$ and $y$ and for all $n\geq 0$. Many such families
are known, and they appear frequently in applications. For example, we
can let $F_n(x)=x^n\!/n!$; the condition
$${(x+y)^n\over n!}=\sum_{k=0}^n\,{x^k\over k!}\;{y^{n-k}\over
(n-k)!}$$
is equivalent to the binomial theorem for integer exponents. Or we can
let $F_n(x)$ be the binomial coefficient ${x\choose n}$; the
corresponding identity
$${x+y\choose n}=\sum_{k=0}^n\,{x\choose k}{y\choose n-k}$$
is commonly called {\it Vandermonde's convolution}.

How special is the convolution condition? {\sl Mathematica\/} 
will readily find all sequences of polynomials that work for,
say, $0\leq n\leq 4$:
\begintt
F[n_,x_]:=Sum[f[n,j]x^j,{j,0,n}]/n!
conv[n_]:=LogicalExpand[Series[F[n,x+y],{x,0,n},{y,0,n}]
     ==Series[Sum[F[k,x]F[n-k,y],{k,0,n}],{x,0,n},{y,0,n}]]
Solve[Table[conv[n],{n,0,4}],
             [Flatten[Table[f[i,j],{i,0,4},{j,0,4}]]]]
\endtt
{\sl Mathematica\/} replies that the $F$'s are either identically zero
or the coefficients of
$F_n(x)=\bigl(f_{n0}+f_{n1}x+f_{n2}x^2+\cdots+f_{nn}x^n\bigr)/n!$ satisfy
$$\eqalign{f_{00}&=1\,,\quad f_{10}=f_{20}=f_{30}=f_{40}=0\,,\cr
\noalign{\smallskip}
f_{22}&=f_{11}^2\,,\quad f_{32}=3f_{11}f_{21}\,,\quad
f_{33}=f_{11}^3\,,\cr
\noalign{\smallskip}
f_{42}&=4f_{11}f_{31}+3f_{21}^2\,,\quad
f_{43}=6f_{11}^2f_{21}\,,\quad f_{44}=f_{11}^4\,.\cr}$$
This allows us to choose $f_{11}$, $f_{21}$, $f_{31}$, and $f_{41}$
freely.

Suppose we weaken the requirements by asking only that the convolution
condition hold when $x=y$. The definition of {\tt conv} then becomes
simply
\begintt
conv[n_]:=LogicalExpand[Series[F[n,2x],{x,0,n}]
            ==Series[Sum[F[k,x]F[n-k,x],{k,0,n}],{x,0,n}]]
\endtt
and we discover that the same solutions occur. In other words, the
weaker requirements imply that the strong requirements are fulfilled
as well.

In fact, it is not difficult to discover a simple rule that
characterizes all ``convolution families.'' {\sl Let
$$F(z)=1+F_1z+F_2z^2+F_3z^3+\cdots$$
 be any power series with\/ $F(0)=1$. Then the polynomials
$$F_n(x)=[z^n]\,F(z)^x$$
form a convolution family. Conversely, every convolution family
arises in this way or is identically zero}. (Here the notation
`$[z^n]$ {\sl expr\/}' stands for what {\sl Mathematica\/} calls {\tt
Coefficient[{\sl expr},z,n]}.)

\bigskip\noindent{\it Proof.}\enspace
Let $f(z)=\ln F(z)=f_1z+f_2z^2\!/2!+f_3z^3\!/3!+\cdots\,$.
It is easy to verify that the coefficient of $z^n$ in $F(z)^x$ is
indeed a polynomial in~$x$ of degree $\leq n$, because 
$F(z)^x=e^{xf(z)}=\exp(xf_1z+xf_2z^2\!/2!+xf_3z^3\!/3!+\cdots\,)$
expands to the power series
$$\sum_{k_1,k_2,k_3,\ldots\geq 0}x^{k_1+k_2+k_3+\cdots}\;
{f_1^{k_1}f_2^{k_2}f_3^{k_3}\ldots\over
1!^{k_1}\,k_1!\,2!^{k_2}\,k_2!\,3!^{k_3}\,k_3!\,\ldots\,}\;
z^{k_1+2k_2+3k_3+\cdots}\,;$$
when $k_1+2k_2+3k_3+\cdots =n$ the coefficient of~$z^n$ is
 a polynomial in~$x$ with terms of
degree $k_1+k_2+k_3+\cdots\leq n$. This construction produces a
convolution family because of the rule for forming coefficients of the
product $F(z)^{x+y}=F(z)^xF(z)^y$.

Conversely, suppose the polynomials $F_n(x)$ form a convolution
family. The condition $F_0(0)=F_0(0)^2$ can hold only if $F_0(x)=0$
or $F_0(x)=1$. In the former case it is easy to prove by induction
that $F_n(x)=0$ for all~$n$. Otherwise, the condition $F_n(0)=2F_n(0)$
for $n>0$ implies that $F_n(0)=0$ for $n>0$. If we equate coefficients
of~$x^k$ on both sides of
$$F_n(2x)=F_n(x)F_0(x)+F_{n-1}(x)F_1(x)+\cdots
+F_1(x)F_{n-1}(x)+F_0(x)F_n(x)\,,$$
we now find that the coefficient $f_{nk}$ of~$x^k$ in $n!\,F_n(x)$ is
forced to have certain values based on the coefficients of
$F_1(x),\ldots,F_{n-1}(x)$, when $k>1$, because $2^kf_{nk}$ occurs on
the left and $2f_{nk}$ on the right. The coefficient $f_{n1}$ can,
however, be chosen freely. Any such choice must make
$F_n(x)=[z^n]\exp(xf_{11}z+xf_{21}z^2\!/2!+xf_{31}z^3\!/3!+\cdots\,)$,
by induction on~$n$.

\bn
{\bf Examples.}\enspace
The first example mentioned above, $F_n(x)=x^n\!/n!$, comes from the
power series $F(z)=e^z$; the second example, $F_n(x)={x\choose n}$,
comes from $F(z)=1+z$. Several other power series are also known to
have simple coefficients when we raise them to the power~$x$.
If $F(z)=1/(1-z)$, for instance, we find
$$[z^n]\,(1-z)^{-x}={-x\choose n}(-1)^n={x+n-1\choose n}\,.$$

It is convenient to use the notations
$$\eqalign{x^{\underline{n}}&=x(x-1)\,\ldots\,(x-n+1)=x!/(x-n)!\cr
\noalign{\smallskip}
x^{\overline{n}}&=x(x+1)\,\ldots\,(x+n-1)=\Gamma(x+n)/\Gamma(x)\cr}$$
for falling factorial powers and rising factorial powers. Since
${x\choose n}=x^{\underline{n}}\!/n!$ and ${x+n-1\choose
n}=x^{\overline{n}}\!/n!$, our last two examples have shown that the
polynomials $x^{\underline{n}}\!/n!$ and $x^{\overline{n}}\!/n!$ form
convolution families, corresponding to $F(z)=1+z$ and $F(z)=1/(1-z)$.
Similarly, the polynomials 
$$F_n(x)={x(x-s)(x-2s)\,\ldots\,\bigl(x-(n-1)\,s)\over n!}$$
form a convolution family corresponding to $(1+sz)^{1/s}$ when $s\neq
0$.

The cases $F(z)=1+z$ and $F(z)=1/(1-z)$ are in fact simply the cases
$t=0$ and $t=1$ of a  general formula for the binomial power
series ${\cal B}_t(z)$, which satisfies
$${\cal B}_t(z)=1+z\,{\cal B}_t(z)^t\,.$$
When $t$ is any real or complex number, exponentiation of this series
is known to yield
$$[z^n]\,{\cal B}_t(z)^x={x+tn\choose n}\,{x\over x+tn}={x(x+tn-1)\ldots
(x+tn-n+1)\over n!}\,;$$
see, for example, 
[Graham et al 1989, section 7.5, example 5],
where a combinatorial proof is given.

The special cases $t=2$ and $t=-1$, 
$$\openup2\jot
\eqalign{{\cal B}_2(z)&={1-\sqrt{\mathstrut 1-4z}\over
2z}=1+z+2z^2+5z^3+14z^4+42z^5+\cdots\;,\cr
{\cal B}_{-1}(z)&={1+\sqrt{\mathstrut 1+4z}\over
2}=1+z-z^2+2z^3-5z^4+14z^5-\cdots\,,\cr}$$
in which the coefficients are the Catalan numbers, arise in numerous
applications. For example, ${\cal B}_2(z)$ is the generating function for
binary trees, and ${\cal B}_1(-z)$ is the reciprocal of ${\cal B}_2(z)$. We
can get identities in trigonometry by noting that $B_2\bigl(({1\over
2}\sin\theta)^2\bigr)={\rm sec}^2(\theta/2)$. Furthermore, if $p$
and~$q$ are probabilities with $p+q=1$, it turns out that
${\cal B}_2(pq)=1/\max(p,q)$. The case $t=1/2$,
$${\cal B}_{1/2}(z)=\left({z+\sqrt{4+z^2}\over 2}\right)^2=1+z+{z^2\over
2}+{z^3\over 2^3}-{z^5\over 2^7}+{2z^7\over 2^{11}}-{5z^9\over
2^{15}}+{14z^{11}\over 2^{19}}-\cdots\,,$$
is another interesting series in which the Catalan numbers can be
seen. The convolution polynomials in this case are the ``central
factorials'' $x(x+{n\over 2}-1)^{\underline{n-1}}/n!$ 
[Riordan 1968, section 6.5], also called Steffensen polynomials
[Roman and Rota 1978, example~6].

The convolution formula corresponding to ${\cal B}_t(z)$,
$${x+y+tn\choose n}\,{x+y\over x+y+tn}=\sum_{k=0}^n\,{x+tk\choose
k}\,{x\over x+tk}\,{y+t(n-k)\choose n-k}\,{y\over y+t(n-k)}$$
is a rather startling generalization of Vandemonde's convolution; it
is an identity for all $x$, $y$, $t$, and~$n$.

The limit of ${\cal B}_t(z/t)^t$ as $t\ra\infty$ is another important
function $T(z)/z$; here
$$T(z)=\sum_{n\geq 1}\,{n^{n-1}\over n!}\,z^n=z+z^2+{3z^3\over
2}+{8z^4\over 3}+{125z^5\over 24}+\cdots$$
is called the {\it tree function\/} because $n^{n-1}$ is the number of
labeled, rooted trees. The tree function satisfies
$$T(z)=ze^{T(z)}\,,$$
and we have the corresponding convolution family
$$[z^n]\left({T(z)\over
z}\right)^x=[z^n]\,e^{xT(z)}={x(x+n)^{n-1}\over n!}\,.$$
The related power series
$$1+zT'(z)={1\over 1-T(z)}=\sum_{n\geq 0}\,{n^nz^n\over n!}=
1+z+2z^2+{9z^3\over 2}+{32z^4\over 3}+{625z^5\over 24}+\cdots$$
defines yet another convolution family of importance: We have
$$[z^n]\,{1\over \bigl(1-T(z)\bigr)^x}={t_n(x)\over n!}\,,$$
where $t_n(x)$ is called the {\it tree polynomial\/} of order~$n$
[Knuth and Pittel 1989].
The coefficients of $t_n(x)=t_{n1}x+t_{n2}x^2+\cdots+t_{nn}x^n$ are
integers with combinatorial significance; namely, $t_{nk}$ is the
number of mappings of an $n$-element set into itself having exactly
$k$~cycles.

A similar but simpler sequence arises from the coefficients of powers
of $e^{ze^z}$:
$$n!\,[z^n]\,e^{xze^z}=\sum_{k=0}^n{n\choose k}k^{n-k}x^k\,.$$
The coefficient of $x^k$ is the number of {\it idempotent\/} mappings
of an $n$-element set into itself, having exactly $k$ cycles
[Harris and Schoenfeld 1967].

If the reader still isn't convinced that convolution families are
worthy of detailed study, well, there's not much hope, although
another example or two might clinch the argument. Consider the power
series
$$e^{e^z-1}=\sum\,{b_nz^n\over n!}=1+{z\over 1!}+{2z^2\over
2!}+{5z^3\over 3!}+{15z^4\over 4!}+{52z^5\over 5!}+\cdots\;;$$
these coefficients $b_n$ are the so-called {\it Bell numbers}, the
number of ways to partition sets of size~$n$ into subsets. For
example, the five partitions that make $b_3=5$ are
$$\{1,2,3\}\,,\quad \{1\}\{2,3\}\,,\quad
\{1,2\}\{3\}\,,\quad\{1,3\}\{2\}\,,\quad \{1\}\{2\}\{3\}\,.$$
The corresponding convolution family is
$$[z^n]\,e^{(e^z-1)x}={{n\brace 0}+{n\brace 1}x+{n\brace
2}x^2+\cdots +{n\brace n}x^n\over n!}\,,$$
where the Stirling number ${n\brace k}$ is the number of partitions
into exactly $k$~subsets.

Need more examples? If the coefficients of $F(z)$ are arbitrary nonnegative
numbers with a finite sum~$S$, then $F(z)/S$ defines a discrete
probability distribution, and the convolution polynomial $F_n(x)$ is
$S^x$~times the probability of obtaining the value~$n$ as the sum of
$x$~independent random variables having that distribution.

\bn
{\bf A derived convolution.}\enspace
Every convolution family $\{F_n(x)\}$ satisfies another general
convolution formula in addition to the one we began with:
$$(x+y)\sum_{k=0}^nk\,F_k(x)\,F_{n-k}(y)=x\,n\,F_n(x+y)\,.$$
For example, if $F_n(x)$ is the convolution family corresponding to
powers of ${\cal B}_t(z)$, this formula says that
$$(x+y)\sum_{k=0}^nk{x+tk\choose k}\;{x\over x+tk}\;{y+t(n-k)\choose
n-k}\;{y\over y+t(n-k)}=xn{x+y+tn\choose n}\;{x+y\over x+y+tn}\,;$$
it looks messy, but it simplifies to another amazing identity in four
parameters,
$$\sum_{k=0}^n{x+t(n-k)\choose n-k}{y+tk\choose k}\;{y\over
y+tk}={x+y+tn\choose n}$$
if we replace $n$ by $n+1$, $k$~by $n+1-k$, and $x$ by $x-t+1$. This
identity  has an interesting history going back to Rothe in 1793
see
[Gould and Kauck\'y 1966].

The alternative convolution formula is proved by differentiating the
basic identity $F(z)^x=\sum_{n\geq 0}F_n(x)z^n$ with respect to~$z$ and
multiplying by~$z$:
$$xzF'(z)\,F(z)^{x-1}=\sum_{n\geq 0}n\,F_n(x)\,z^n\,.$$
Now $\sum_{k=0}^n k\,F_k(x)\,F_{n-k}(y)$ is the coefficient of~$z^n$ in
$xzF'(z)\,F(z)^{x+y-1}$, while $nF_n(x+y)$ is the coefficient of~$z^n$
in $(x+y)zF'(z)\,F(z)^{x+y-1}$.
Q.E.D.

\bn
{\bf Convolution and composition.}\enspace
Once upon a time I was trying to remember the form of a general
convolution family, so I~gave {\sl Mathematica\/} the following command:
\begintt
Simplify[Series[(1+Sum[A[k]z^k,{k,4}])^x,{z,0,4}]]
\endtt
The result was a surprise. Instead of presenting the coefficient
of~$z^n$ as a polynomial in~$x$, {\sl Mathematica\/} chose another
form: The coefficient of~$z^2$, for example, was ${1\over
2}A_1^2x(x-1)+A_2x$. In the notation of
 falling factorial powers, {\sl Mathematica}'s response took the form
$$\eqalign{&\textstyle{1+A_1xz+\bigl({1\over
2}A^2_1x^{\underline{2}}+A_2x\bigr)\,z^2+\bigl({1\over
6}A_1^3x^{\underline{3}} +A_1A_2x^{\underline{2}}+A_3x\bigr)\,z^3}\cr
\noalign{\smallskip}
&\textstyle{\qquad\null+\bigl({1\over 24}A_1^4x^{\underline{4}}+{1\over
2}A_1^2A_2x^{\underline{3}}+\bigl(A_1A_3+{1\over
2}A_2^2\bigr)\,x^{\underline{2}}+ A_4x\bigr)\,z^4+O(z)^5}\,.\cr}$$
I wasn't prepared to work with factorial powers, so I tried another
tack:
\begintt
Simplify[Series[Exp[Sum[a[k]z^k,{k,4}]x],{z,0,4}]]
\endtt
This time I got ordinary polynomials in $x$, but---lo and
behold---they were
$$\eqalign{&\textstyle{%
1+a_1xz+\bigl({1\over 2}a_1^2x^2+a_2x\bigr)\,z^2+\bigl({1\over
6}a_1^3x^3+a_1a_2x^2+a_3x)\,z^3}\cr
\noalign{\smallskip}
&\textstyle{\qquad\null+\bigl({1\over 24}a_1^4x^4+{1\over
2}a_1^2a_2x^3+\bigl(a_1a_3+{1\over 2}a_2^2\bigr)x^2+a_4x\bigr)\,z^4
+O(z)^5}\,.\cr}$$
The result was exactly the same as before, but with $a$'s in place
of~$A$'s, and with normal powers in place of the factorials!

So I learned a curious phenomenon: {\sl If we take any convolution
family and replace each power\/ $x^k$ by\/~$x^{\underline{k}}$, we get
another convolution family}. (By the way, the replacement can be done
in {\sl Mathematica\/} by saying
\begintt
Expand[F[n,x]]/.Power[x,k_]->k!Binomial[x,k];
\endtt
expansion is necessary in case $F_n(x)$ has been factored.)

The proof was not difficult to find, once I psyched out how {\sl
Mathematica\/} might have come up with its factorial-based formula:
We have
$$e^{xf(z)}=1+f(z)\,x+{f(z)^2\over 2!}\,x^2+{f(z)^3\over 3!}\,x^3
+\cdots\;,$$
and furthermore
$$\bigl(1+f(z)\bigr)^x=1+f(z)\,x+{f(z)^2\over
2!}\;x^{\underline{2}}+{f(z)^3\over 3!}\;x^{\underline{3}}+\cdots\;.$$
Therefore if we start with the convolution family $F_n(x)$
corresponding to $F(z)=e^{f(z)}$, and replace each~$x^k$
by~$x^{\underline{k}}$, we get the convolution family corresponding to
$1+f(z)=1+\ln F(z)$.

A similar derivation shows that if we replace $x^k$ by the {\it
rising\/} factorial power
$x^{\overline{k}}$
instead, we get the convolution family corresponding to
$1/\bigl(1-f(z)\bigr)=1/\bigl(1-\ln F(z)\bigr)$. 
In particular, if we begin with the family $F_n(x)=x(x+n)^{n-1}\!/n!$
corresponding to $T(z)/z=e^{T(z)}$, and if we replace~$x^k$
by~$x^{\overline{k}}$ to get
$${1\over n!}\,\sum_{k=0}^{n-1}{n-1\choose k}\,x^{\overline{k+1}}\,
n^{n-1-k}\,,$$
this must be $[z^n]\,\bigl(1-T(z)\bigr)^{-x}=t_n(x)/n!$, the tree
polynomial. 

Indeed, we can replace each $x^k$ by $k!\,G_k(x)$, where $\{G_k(x)\}$
is any convolution family whatever!  The previous examples,
$x^{\underline{k}}$ and~$x^{\overline{k}}$, are merely the special
cases $k!\,{x\choose k}$ and $k!\,{x+k-1\choose k}$ corresponding to
two of the simplest and most basic families we have considered. In
general we get
$$1+f(z)\,G_1(x)+{f(z)^2\over 2!}\;2!\,G_2(x)+{f(z)^3\over
3!}\,3!\,G_3(x)+\cdots\;,$$ 
which is none other than $G\bigl(f(z)\bigr)^x=G\bigl(\ln
F(z)\bigr)^x$. 

For example, $G_k(x)={x+2k\choose k}\;{x\over
x+2k}=x(x+2k-1)^{\underline{k-1}}/k!$ is the family corresponding to
${\cal B}_2(z)$. If we know the family $F_n(x)$ corresponding to
$e^{f(z)}$ we can replace~$x^k$ by $x(x+2k-1)^{\underline{k-1}}$,
thereby obtaining the family that corresponds to
${\cal B}_2\bigl(f(z)\bigr)=\bigl(1+\sqrt{1-4f(z)}\,\bigr)/2f(z)$.

\bn
{\bf Convolution matrices.}\enspace
I knew that such remarkable facts must have been discovered before,
although they were new to me at the time. And indeed, it was not
difficult to find them in books, once I~knew what to look for.
(Special cases of general theorems are not always easy to recognize,
because any particular formula is a special case of infinitely many
generalizations, almost all of which are false.)

In the special case that each polynomial $F_n(x)$ has degree
exactly~$n$, i.e., when $f_1\neq 0$, the polynomials $n!\,F_n(x)$ are
said to be of {\it binomial type\/}
[Mullin and Rota 1970].
An extensive theory of such polynomial sequences has been developed 
[Rota et al 1973] [Garsia 1973] [Roman and Rota 1978],
based on the theory of linear operators, and the reader will find it
quite interesting to compare the instructive treatment in those papers
to the related but rather different directions explored in the present
work. A~comprehensive exposition of the operator approach appears in
[Roman 1984]. Actually, Steffensen had defined a concept called {\it
poweroids}, many years earlier [Steffensen 1941], and poweroids are
almost exactly the same as sequences of binomial type; but Steffensen
apparently did not realize that his poweroids satisfy the convolution
property, which we can readily deduce (with hindsight) from
equations~(6) and~(7) of his paper.

Eri Jabotinsky introduced a nice way to understand the phenomena of
convolution polynomials, by considering the infinite matrix of
coefficients~$f_{nk}$ [Jabotinsky 1947].
Let us
recapitulate the notation that was introduced informally above:
$$\eqalign{e^{xf(z)}&=F(z)^x=1+F_1(x)\,z+F_2(x)\,z^2+\cdots\,;\cr
\noalign{\smallskip}
F_n(x)&=(f_{n1}x+f_{n2}x^2+\cdots+f_{nn}x^n)/n!\,;\cr
\noalign{\smallskip}
f(z)&=f_1z+f_2z^2\!/2!+f_3z^3\!/3!+\cdots\;.\cr}$$
Then Jabotinsky's matrix $F=(f_{nk})$ is a lower triangular matrix
containing the coefficients of $n!\,F_n(x)$ in the $n$th row. The
first few rows are
$$\vcenter{\halign{$#$\hfil\quad&\hfil$#$\hfil\quad&\hfil$#$\hfil\quad
&\hfil$#$\hfil\cr
f_1\cr
\noalign{\smallskip}
f_2&\phantom{1}f_1^2\cr
\noalign{\smallskip}
f_3&3f_1f_2&\phantom{1}f_1^3\cr
\noalign{\smallskip}
f_4&4f_1f_3+3f_2^2&6f_1^2f_2&f_1^4\,,\cr}}$$
as we saw earlier. In general,
$$f_{nk}=\sum\,{n!\over
1!^{k_1}\,k_1!\,2!^{k_2}\,k_2!\,3!^{k_3}\,k_3!\,\ldots}\; 
f_1^{k_1}f_2^{k_2}f_3^{k_3}\ldots\;,$$
summed over all $k_1,k_2,k_3,\ldots\geq 0$ with
$$k_1+k_2+k_3+\cdots =k\,,\qquad k_1+2k_2+3k_3+\cdots =n\,.$$
(The summation is over all partitions of the integer~$n$ into
$k$~parts, where $k_j$ of the parts are equal to~$j$.) We will call
such an array a {\it convolution matrix}.

If each original coefficient $f_j$ is an integer, all entries of the
corresponding convolution matrix will be integers, because the
complicated quotient of factorials in the sum is an integer---it is
the number of ways to partition a set of $n$~elements into $k$~subsets
with exactly $k_j$ of the subsets having size~$j$. Given the first
column we can compute the other columns from left to right and from
top to bottom by using the recurrence
$$f_{nk}=\sum_{j=1}^{n-k+1}\,{n-1\choose j-1}f_j\,f_{(n-j)(k-1)}\,.$$
This recurrence is based on set partitions on which the element~$n$
occurs in a subset of size~$j$: There are ${n-1\choose j-1}$ ways to
choose the other $j-1$ elements of the subset, and the factor
$f_{(n-j)(k-1)}$ corresponds to partitioning the remaining $n-j$
elements into $k-1$ parts.

For example, if each $f_j=1$, the convolution matrix begins
$$\vcenter{\halign{\hfil$#$\quad&\hfil$#$\quad&\hfil$#$\quad
&\hfil$#$\quad&\hfil$#$\quad&$#$\hfil\cr
1\cr
1&1\cr
1&3&1\cr
1&7&6&1&&.\cr
1&15&25&10&1\cr}}$$
These are the numbers ${n\brace k}$ that {\sl Mathematica\/} calls
{\tt StirlingS2[n,k]}; they arose in our example of Bell numbers when
$f(z)=e^z-1$. Similarly, if each $f_j=(j-1)!$, the first five rows are
$$\vcenter{\halign{\hfil$#$\quad&\hfil$#$\quad&\hfil$#$\quad
&\hfil$#$\quad&\hfil$#$\quad&\hfil$#$\cr
1\cr
1&1\cr
2&3&1\cr
6&11&6&1&&;\cr
24&50&35&10&1\cr}}$$
{\sl Mathematica\/} calls these numbers {\tt
(-1){\char'136}(n-k)StirlingS1[n,k]}. In this case
$f(z)=\ln\bigl(1/(1-z)\bigr)$, and $F_n(z)={x+n-1\choose n}$.
The signed numbers {\tt StirlingS1[nk]}, 
$$\vcenter{\halign{\hfil$#$\quad&\hfil$#$\quad&\hfil$#$\quad
&\hfil$#$\quad&\hfil$#$\quad&\hfil$#$\cr
1\cr
{-}1&1\cr
2&{-}3&1\cr
{-}6&11&{-}6&1\cr
24&{-}50&35&{-}10&1\cr}}$$
correspond to $f(z)=\ln(1+z)$ and $F_n(z)={x\choose n}$. In general if
we replace~$z$ by~$\alpha z$ and $x$ by~$\beta x$, the effect is to
multiply row~$n$ of the matrix by~$\alpha^n$ and to multiply
column~$k$ by~$\beta^k$. Thus when $\beta=\alpha^{-1}$, the net effect
is to multiply $f_{nk}$ by~$\alpha^{n-k}$. Transforming the signs by a
factor $(-1)^{n-k}$ corresponds to changing $F(z)$ to $1/F(-z)$ and
$f(z)$ to $-f(-z)$. Therefore the matrix that begins
$$\vcenter{\halign{\hfil$#$\quad&\hfil$#$\quad&\hfil$#$\quad
&\hfil$#$\quad&\hfil$#$\quad&\hfil$#$\cr
1\cr
{-}1&1\cr
1&{-}3&1\cr
{-}1&7&{-}6&1\cr
1&{-}15&25&{-}10&1\cr}}$$
corresponds to $f(z)=1-e^{-z}$.

Let's look briefly at some of our other examples in matrix form. When
$F(z)={\cal B}_t(z)$, we have $f_j=(tj-1)^{\underline{j-1}}$, which is an
integer when $t$ is an integer. In particular, the Catalan case $t=2$
produces a matrix that begins
$$\vcenter{\halign{\hfil$#$\quad&\hfil$#$\quad&\hfil$#$\quad
&\hfil$#$\quad&\hfil$#$\quad&\hfil$#$\cr
1\cr
3&1\cr
20&9&1\cr
210&107&18&1\cr
3024&1650&335&30&1&.\cr}}$$
When $t=1/2$, we can remain in an all-integer realm by replacing $z$
by~$2z$ and~$x$ by~$x/2$. Then $f_j=0$ when $j$ is even, while
$f_{2j+1}=(-1)^j(2j-1)!!^2$:
$$\vcenter{\halign{\hfil$#$\quad&\hfil$#$\quad&\hfil$#$\quad
&\hfil$#$\quad&\hfil$#$\quad&\hfil$#$\cr
1\cr
0&1\cr
{-}1&0&1\cr
0&{-}4&0&1&&;\cr
9&0&{-}10&0&1\cr}}$$
If we now replace $z$ by $iz$ and $x$ by $x/i$ to eliminate the minus
signs, we find that 
$f(z)=$ arcsin~$z$, because
$\ln\bigl(iz+\sqrt{1-z^2}\,\bigr)=i\theta$ when
$z=\sin\theta$. Thus we can deduce a closed
form for the coefficients of $e^{x\arcsin
z}={\cal B}_{1/2}(2iz)^{x/(2i)}$:
$$\eqalign{%
n!\,[z^n]\,e^{x\arcsin z}&=(2i)^{n-1}x\left({x\over 2i}+{n\over
2}-1\right)\,\ldots \left({x\over 2i}-{n\over 2}+1\right)\cr
\noalign{\smallskip}
&=\cases{x^2(x^2+2^2)\,\ldots\,\bigl(x^2+(n-2)^2\bigr)\,,&$n$ even;\cr
\noalign{\smallskip}
x(x^2+1^2)(x^2+3^2)\,\ldots\,\bigl(x^2+\bigl((n-2)^2\bigr)\bigr)\,,&$n$
odd.\cr}\cr}$$
This remarkable formula is equivalent to the theorem of [Gomes
Teixeira 1896].

If $f_j=2^{1-j}$ when $j$ is odd but $f_j=0$ when $j$ is even, we get
the convolution matrix corresponding to $e^{2x\sinh(z/2)}$:
$$\vcenter{\halign{$\hfil#\hfil$\quad
&$\hfil#\hfil$\quad
&$\hfil#\hfil$\quad
&\hfil$#$\quad
&$\hfil#\hfil$\quad
&\hfil$#$\quad
&$\hfil#\hfil$\quad
&$\hfil#\hfil$\quad
&$\hfil#\hfil$\cr
1\cr
0&1\cr
{1\over 4}&0&1\cr
0&1&0&1\cr
{1\over 16}&0&{5\over 2}&0&1\cr
0&1&0&5&0&1\cr
{1\over 64}&0&{91\over 16}&0&{35\over 4}&0&1\cr
0&1&0&21&0&14&0&1&.\cr}}$$
Again we could stay in an all-integer realm if we replaced $z$ by~$2z$
and~$x$ by~$x/2$; but the surprising thing in this case is that the
entries in even-numbered rows and columns are all integers {\it
before\/} we make any such replacement. The reason is that the entries
satisfy $f_{nk}=k^2f_{(n-2)k}/4+f_{(n-2)(k-2)}$. (See [Riordan 1968,
pages 213--217], where the notation $T(n,k)$ is used for these
``central factorial numbers''~$f_{nk}$.)

We can complete our listing of noteworthy examples by setting
$f_j=\sum_{k=1}^n\,n^{n-k-1}n^{\underline{k}}$; then we get the
coefficients of the tree polynomials:
$$\vcenter{\halign{\hfil$#$\quad&\hfil$#$\quad&\hfil$#$\quad
&\hfil$#$\quad&\hfil$#$\quad&\hfil$#$\cr
1\cr
3&1\cr
17&9&1\cr
142&95&18&1\cr
1569&1220&305&30&1&.\cr}}$$
The sum of the entries in row $n$ is $n^n$.

\bn
{\bf Composition and iteration.}\enspace
Jabotinski's main reason for defining things as he did was
his observation that {\sl the product of convolution matrices is a
convolution matrix}.
Indeed, if $F$ and $G$ are the convolution matrices corresponding to
the functions $e^{xf(z)}$ and $e^{xg(z)}$ we have the vector/matrix
identities 
$$\eqalign{e^{xf(z)}-1
&=(z,z^2\!/2!,z^3\!/3!,\ldots\,)\,F\,(x,x^2,x^3,\ldots\,)^{\rm T}\cr
\noalign{\smallskip}
e^{xg(z)}-1
&=(z,z^2\!/2!,z^3\!/3!,\ldots\,)\,G\,(x,x^2,x^3,\ldots\,)^{\rm T}\cr}$$ 
If we now replace $x^k$ in $e^{xf(z)}$ by $k!\,G_k(x)$, as in our
earlier discussion, we get
$$\eqalign{(z,z^2\!/2!,z^3\!/3!,\ldots\,)\,F\,&\bigl(G_1(x),2!\,G(x),
3!\,G_3(x),\ldots\,\bigr)^{\rm T}\cr
\noalign{\smallskip}
=&(z,z^2\!/2!,z^3\!/3!,\ldots\,)\,FG\,(x,x^2,x^3,\ldots\,)^{\rm T}\cr
\noalign{\smallskip}
=&\bigl(f(z),f(z)^2\!/2!,f(z)^3\!/3!,
\ldots\,\bigr)\,G\,(x,x^2,x^3,\ldots\,)^{\rm T}\cr
\noalign{\smallskip}
=&e^{xg(f(z))}-1\,.\cr}$$
Multiplication of convolution matrices corresponds to composition of
the functions in the exponent.

Why did the function corresponding to~$FG$ turn out to be
$g\bigl(f(z)\bigr)$ instead of $f\bigl(g(z)\bigr)$? Jabotinsky, in
fact, defined his matrices as the transposes of those given here. The
rows of his (upper triangular) matrices were the power series $f(z)^k$,
while the columns were the polynomials $F_n(x)=[z^n]\,e^{xf(z)}$; with
those conventions the product of his matrices $F^{\rm T}G^{\rm T}$
 corresponded to
$f\bigl(g(z)\bigr)$. (In fact, he defined a considerably more general
representation, in which the matrix~$F$ could be $U^{-1}FU$ for any
nonsingular matrix~$U$.)
However, when our interest is focussed on the polynomials
$n!\,F_n(x)$, as when we study Stirling numbers or tree polynomials or
the Stirling polynomials to be discussed below, it is more natural to
work with lower triangular matrices and to insert factorial
coefficients, as Comtet did [Comtet 1970, section 3.7]. The two
conventions are isomorphic. Without the factorials, convolution
matrices are sometimes called {\it renewal arrays\/} [Rogers 1978].
We would get a non-reversed order if we
had been accustomed to using postfix notation $(z)f$ for functions, as
we do for operations such as squaring or taking transposes or
factorials; then $g\bigl(f(z)\bigr)$ would be $\bigl((z)f\bigr)g$. 

Recall that the Stirling numbers ${n\brace k}$ correspond to
$f(z)=e^z-1$, and the Stirling numbers ${n\brack k}$ correspond to
$g(z)=\ln\bigl(1/(1-z)\bigr)$. Therefore if we multiply Stirling's
triangles we get the convolution matrix
$$\vcenter{\halign{\hfil$#\;$&\hfil$#$\quad&\hfil$#$\quad&\hfil$#$\quad
&\hfil$#$\quad&\hfil$#$\quad&\hfil$#$\cr
&1\cr
&2&1\cr
FG=&6&6&1&&&,\cr
&26&36&12&1\cr
&150&250&120&20&1\cr}}$$
which corresponds to $g\bigl(f(z)\bigr)=\ln\bigl(1/(2-e^z)\bigr)$.
Voila!
These convolution polynomials represent the coefficients of
$(2-e^z)^{-x}$. [Cayley 1859] showed that $(2-e^z)^{-1}$ is the
exponential generating function for the sequence
$1,3,13,75,541,\ldots\,$, which counts {\it preferential
arrangements\/} of $n$~objects, i.e., different outcomes of sorting
when equality is possible as well as inequality. The coefficient
$(fg)_{nk}$ is the number of preferential arrangements in which the
``current minimum'' changes $k$~times when we examine the elements one
by one in some fixed order. (See [Graham et al 1989, exercise 7.44].)

Similarly, the reverse matrix product  yields the so-called Lah numbers
[Lah 1955],
$$\vcenter{\halign{\hfil$#\;$&\hfil$#$\quad&\hfil$#$\quad&\hfil$#$\quad
&\hfil$#$\quad&\hfil$#$\quad&\hfil$#$\cr
&1\cr
&2&1\cr
GF=&6&6&1&&&;\cr
&24&36&12&1\cr
&120&240&120&20&1\cr}}$$
here $f_j=j!$ and the rows represent the coefficients of
$\exp\bigl(xf\bigl(g(z)\bigr)\bigr)=\exp(xz+xz^2+xz^3+\cdots\,)$.
Indeed, the convolution polynomials in this case are the generalized
Laguerre polynomials $L_n^{(-1)}(-x)$, which {\sl Mathematica\/} calls
{\tt LaguerreL[n,-1,-x]}. These polynomials can also be expressed as
$L_n(-x)-L_{n-1}(-x)$; or as {\tt LaguerreL[n,-x]-LaguerreL[n-1,-x]} if
we say
\begintt
Unprotect[LaguerreL]; LaguerreL[-1,x_]:=0; Protect[LaguerreL]
\endtt
first.
The row sums $1,3,13,73,501,\ldots$ of~$GF$ 
enumerate ``sets of lists'' 
[Motzkin 1971]; 
the coefficients are
$(GF)_{nk}=n!\,[z^n]\,f\bigl(g(z)\bigr)^k\!/k!={n\choose k}{n-1\choose
k-1}(n-k)!$ [Riordan 1968, exercise 5.7].

Since convolution matrices are closed under multiplication, they are
also closed under exponentiation, i.e., under taking of powers. The
$q$th power~$F^q$ of a convolution matrix then corresponds to $q$-fold
iteration of the function $\ln F=f$. Let us denote $f\bigl(f(z)\bigr)$
by $f^{[2]}(z)$; in general, the $q$th iterate $f^{[q]}(z)$ is defined
to be $f\bigl(f^{[q-1]}(z)\bigr)$, where $f^{[0]}(z)=z$. This is {\sl
Mathematica}'s {\tt Nest[f,z,q]}.

The $q$th iterate can be obtained by doing $O(\log q)$ matrix
multiplications, but in the interesting case $f'(0)=f_1=1$ we can also
compute the coefficients of $f^{[q]}(z)$ by using formulas in which
$q$ is simply a numerical parameter. Namely, as suggested by
[Jabotinsky 1947], we can express the matrix power~$F^q$ as
$$\bigl(I+(F-I)\bigr)^q=I+{q\choose 1}(F-I)+{q\choose
2}(F-I)^2+{q\choose 3}(F-I)^3+\cdots\;.$$
This infinite series converges, because the entry in row $n$ and
column~$k$ of $(F-I)^j$ is zero for all $j>n-k$. When $q$ is any
positive integer, the result defined in this way is a convolution
matrix. Furthermore, the matrix entries are all polynomials in~$q$.
Therefore the matrix obtained by this infinite series is a convolution
matrix for all values of~$q$.

Another formula for the entries of $F^q$ was presented in 
[Jabotinsky 1963]. 
Let $f_{nk}^{(q)}$ be the element in row~$n$ and column~$k$; then
$$\openup3\jot\eqalign{f_{nk}^{(q)}
&=\sum_{l=0}^m{q\choose l}(F-I)_{nk}^l\cr
&=\sum_{j=0}^mf_{nk}^{(j)}\sum_{l=j}^m{q\choose l}{l\choose
j}(-1)^{l-j}\cr
&=\sum_{j=0}^mf_{nk}^{(j)}{q\choose j}\sum_{l=j}^m{q-j\choose
l-j}(-1)^{l-j}\cr
&=\sum_{j=0}^mf_{nk}^{(j)}{q\choose j}{q-j-1\choose
m-j}(-1)^{m-j}\,,\cr}$$
for any $m\geq n-k$.
Indeed, we have $p(q)=\sum_{j=0}^mp(j){q\choose j}{q-j-1\choose
m-j}(-1)^{m-j}$ whenever $p$ is a polynomial of degree $\leq m$; this
is a special case of Lagrange interpolation.

It is interesting to set $q=1/2$ and compute
convolution square roots of the Stirling number matrices. We have
$$\eqalign{\pmatrix{1\cr 1/2&1\cr 1/8&3/2&1\cr 0&5/4&3&1\cr
1/32&5/8&5&5&1\cr}^2
&=\pmatrix{1\cr 1&1\cr 1&3&1\cr 1&7&6&1\cr 1&15&25&10&1\cr}\,;\cr
\noalign{\medskip}
\pmatrix{1\cr 1/2&1\cr 5/8&3/2&1\cr 5/4&13/4&3&1\cr
109/32&75/8&10&5&1\cr}^2 
&=\pmatrix{1\cr 1&1\cr 2&3&1\cr 6&11&6&1\cr 24&50&35&10&1\cr}\,.\cr}$$
The function $z+z^2\!/4+z^3\!/48+z^5\!/3840-7z^6\!/92160+\cdots$ 
therefore lies
``halfway'' beween~$z$ and $e^z-1=z+z^2\!/2!+z^3\!/3!+\cdots\;$, and the 
function
$z+z^2\!/4+5z^3\!/48+5z^4\!/96+109z^5\!/3840+497z^6\!/30720+\cdots$
lies halfway between $z$ and $\ln 1/(1-z)=z+z^2\!/2+z^3\!/3+\cdots\;$.
These half-iterates are unfamiliar functions; but it is not difficult
to prove that $z/(1-z/2)=z+z^2\!/2+z^3\!/4+\cdots$ is halfway
between~$z$ and $z/(1-z)=z+z^2+z^3+\cdots\,$. In general when
$f(z)=z/(1-cz^k)^{1/k}$ we have $f^{[q]}(z)=z/(1-qcz^k)^{1/k}$. 
 
It seems natural to conjecture that the coefficients of $f^{[q]}(z)$
are positive for $q>0$ when $f(z)=\ln 1/(1-z)$; but this conjecture
turns out to be false, because {\sl Mathematica\/} reports that
$[z^8]\,f^{[q]}(z)=-11q/241920+O(q^2)$. Is there a simple necessary
and sufficient condition on~$f$ that characterizes when all
coefficients of~$f^{[q]}$ are nonnegative for nonnegative~$q$? This
will happen if and only if the entries in the first column of
$$\ln F=(F-I)-{\textstyle{1\over 2}}(F-I)^2+{\textstyle{1\over 3}}
(F-I)^3-\cdots$$ 
are nonnegative. (See [Kuczma 1968] for iteration theory and an
extensive bibliography.)

\bn
{\bf Reversion.}\enspace
The case $q=-1$ of iteration is often called reversion of series,
although {\sl Mathematica\/} uses the more proper name 
{\tt InverseSeries}. Given $f(z)=f_1z+f_2z^2\!/2!+\cdots\,$, we seek
$g(z)=f^{[-1]}(z)$ such that $g\bigl(f(z)\bigr)=z$. This is clearly
equivalent to finding the first column of the
inverse of the convolution matrix. 

The inverse does not exist when $f_1=0$, because the diagonal of~$F$
is zero in that case. Otherwise we  can assume that $f_1=1$, because
$f_1g\bigl(f(z/f_1)\bigr)=z$ when $g$ reverts the power series
$f(z/f_1)$.

When $f_1=1$ we can obtain the inverse by setting $q=-1$ in our
general formula for iteration. But Lagrange's celebrated inversion
theorem for power series tells us that there is another, more
informative, way to compute the function $g=f^{[-1]}$.
Let us set $\widehat{F}(z)=f(z)/z=1+f_2z/2!+f_3z^2\!/3!+\cdots\,$. Then
Lagrange's theorem states that the elements of the matrix $G=F^{-1}$
are
$$g_{nk}={(n-1)!\over (k-1)!}\,\widehat{F}_{n-k}(-n)\,,$$
where $\widehat{F}_n(x)$ denotes the convolution family corresponding to
$\widehat{F}(z)$.

There is a surprisingly simple way to prove Lagrange's theorem, using
our knowledge of convolution families. Note first that
$$f_{nk}=n!\,[z^nx^k]\,e^{xf(z)}={n!\over k!}\;[z^n]\;f(z)^k={n!\over
k!}\;[z^{n-k}]\;\widehat{F}(z)^k\,;$$
therefore
$$f_{nk}={n!\over k!}\;\widehat{F}_{n-k}(k)\,.$$
Now we need only verify that the matrix product $GF$ is the identity,
by computing its element in row~$n$ and column~$m$:
$$\sum_{k=m}^ng_{nk}f_{km}=\sum_{k=m}^n\;{(n-1)!\over
(k-1)!}\;\widehat{F}_{n-k}\,(-n)\;{k!\over m!}\;\widehat{F}_{k-m}(m)\,.$$
When $m=n$ the sum is obviously 1. When $m=n-p$ for $p>0$ it is
$(n-1)!/(n-p)!$ times 
$$\openup3\jot
\eqalign{\sum_{k=n-p}^nk\,\widehat{F}_{n-k}(-n)\,\widehat{F}_{k-n+p}(n-p)
&=\sum_{k=0}^p\,(n-k)\,\widehat{F}_k(-n)\,\widehat{F}_{p-k}(n-p)\cr
&=n\,\sum_{k=0}^p\,\widehat{F}_k(-n)\,\widehat{F}_{p-k}(n-p)-\sum_{k=0}^p\,k\,
\widehat{F}_k(-n)\,\widehat{F}_{p-k}(n-p)\cr
&=n\,\widehat{F}_p(-p)-n\,\widehat{F}_p(-p)=0\cr}$$
by the original convolution formula and the one we derived from it.
The proof is complete.

\bn
{\bf Extending the matrix.}\enspace
The simple formula for $f_{nk}$ that we used to prove Lagrange's
theorem when $f_1=1$ can be written in another suggestive form, if we
replace $k$ by $n-k$:
$$f_{n(n-k)}=n^{\underline{k}}\,\widehat{F}_k(n-k)\,.$$
For every fixed $k$, this is a polynomial in~$n$, of degree $\leq 2k$.
Therefore we can define the quantity $f_{y(y-k)}$ for all real or
complex~$y$ to be $y^{\underline{k}}\,\widehat{F}_k(y-k)$; and in
particular we can define $f_{nk}$ in this manner for all integers $n$
and~$k$, letting $f_{nk}=0$ when $k>n$. For example, in the case of
Stirling numbers this analysis establishes the well-known fact
 that ${y\brace y-k}$ and ${y\brack
y-k}$ are polynomials in~$y$ of degree~$2k$, and that these
polynomials are multiples of
$y^{\underline{k+1}}=y(y-1)\,\ldots\,(y-k)$ when $k>0$.

The two flavors of Stirling numbers are related in two important ways.
First,  their matrices are inverse to each other if we
attach the signs $(-1)^{n-k}$ to the elements in one matrix:
$$\sum_{k=0}^n\,{n\brace k}{k\brack m}(-1)^{n-k}=\sum_{k=0}^m
{n\brack k}{k\brace m}(-1)^{n-k}=\delta_{mn}\,.$$
This follows since the numbers ${n\brace k}$ correspond to
$f(z)=e^z-1$ and the numbers ${n\brack k}(-1)^{n-k}$  correspond to
$g(z)=\ln(1+z)$, as mentioned earlier, and we have $g\bigl(f(z)\bigr)=z$.

The other important relationship beween ${n\brace k}$ and ${n\brack
k}$ is the striking identity
$${n\brace k}={-k\brack -n}\,,$$
which holds for all integers $n$ and $k$ when we use the polynomial
extension method. We can prove in fact, that the analogous relation
$$f_{nk}=(-1)^{k-n}g_{(-k)(-n)}$$
holds in the extended matrices $F$ and $G$ that correspond to {\it
any\/} pair of inverse functions $g\bigl(f(z)\bigr)=z$, when
$f'(0)=1$. For we have
$$(-1)^{n-k}g_{(-k)(-n)}=(-1)^{n-k}(-k-1)(-k-2)\,
\ldots\,(-n)\,\widehat{F}_{n-k}(k)
={n!\over k!}\;\widehat{F}_{n-k}(k)=f_{nk}$$
in the formulas above. (The interesting history of the identity
${n\brace k}={-k\brack -n}$ is traced in [Knuth 1992]. The fact that
the analogous
 formula holds in any convolution matrix was pointed out by Ira
Gessel after he had read a draft of that paper. See also 
[Jabotinski 1953]; [Carlitz 1978];
[Roman and Rota 1978, section~10].)

Suppose we denote the Lah numbers ${n\choose k}{n-1\choose k-1}(n-k)!$
by $\left\vert{n\atop k}\right\vert$. The extended matrix in that case
has a pleasantly symmetrical property
$$\left\vert{n\atop k}\right\vert=\left\vert{-k\atop
-n}\right\vert\,,$$
because the corresponding function $f(z)=z/(1-z)$ satisfies
$f\bigl(-f(-z)\bigr)=z$. (Compare
[Mullin and Rota 1969, section~9].)
Near the origin $n=k=0$, the nonzero entries
look like this:
$$\vcenter{\halign{
\hfil#\quad
&\hfil#\quad
&\hfil#\quad
&\hfil#\quad
&\hfil#\quad
&\hfil#\quad
&\hfil#\quad
&\hfil#\quad
&\hfil#\quad
&\hfil#\quad
&\hfil#\cr
$\ldots$&1\cr
$\ldots$&12&1\cr
$\ldots$&36&6&1\cr
$\ldots$&24&6&2&1\cr
&&&&&1\cr
&&&&&&1\cr
&&&&&&2&1\cr
&&&&&&6&6&1\cr
&&&&&&\kern-.5em24&36&12&1\cr
&&&&&&$\vdots$&$\vdots$&$\vdots$\cr
}}$$

\bn
{\bf Still more convolutions.}\enspace
Our proof of Lagrange's theorem yields yet another corollary. Suppose
$g\bigl(f(z)\bigr)=z$ and $f'(0)=1$, and let $\widehat{F}(z)=f(z)/z$,
$\widehat{G}(z)=g(z)/z$. Then the equation
$$g_{nk}={n!\over k!}\;\widehat{G}_{n-k}(k)={(n-1)!\over
(k-1)!}\;\widehat{F}_{n-k}(-n)$$ 
tell us, after replacing $n$ by $n+k$, that the identity
$${n+k\over k}\;\widehat{G}_n(k)=\widehat{F}_n(-n-k)$$
holds for all positive integers~$k$. Thus the polynomials
$\widehat{G}_n(x)$ and $\widehat{F}_n(x)$ must be related by the formula
$$(x+n)\,\widehat{G}_n(x)=x\widehat{F}_n(-x-n)\,.$$

Now $\widehat{F}_n(x)$ is an arbitrary convolution family, and
$\widehat{F}_n(-x)$ is another. We can conclude that {\sl if\/
$\{F_n(x)\}$ is any convolution family, then so is the set of
polynomials\/} $\{xF_n(x+n)/(x+n)\}$.
Indeed, if $F_n(x)$ corresponds to the coefficients of $F(z)^x$, our
argument proves that the coefficients of $G(z)^x$ are
$x\,F_n(x+n)/(x+n)$, where $zG(z)$ is the inverse of the power series
$z/F(z)$:
$$G(z)=F\bigl(zG(z)\bigr)\,,\qquad G\bigl(z/F(z)\bigr)=F(z)\,.$$
The case $F(z)=1+z$ and $G(z)=1/(1-z)$ provides a simple example,
where we know that $F_n(x)={x\choose n}$ and $G_n(x)={x+n-1\choose
n}=xF_n(x+n)/(x+n)$.

A more interesting example arises when
$F(z)=ze^z\!/(e^z-1)=z+z/(e^z-1)
=1+z/2+B_2z^2\!/2!+B_4z^4\!/4!+\cdots\,$; then $F(-z)$ is the
exponential generating function for the Bernoulli numbers. The
convolution family for $F(z)^x$ is $F_n(x)=x\sigma_n(x)$, where
$\sigma_n(x)$ is called a {\it Stirling polynomial}. (Actually
$\sigma_0(x)=1/x$, but $\sigma_n(x)$ is a genuine polynomial when
$n\geq 1$.) The function~$G$ such that $G\bigl(z/F(z)\bigr)=F(z)$ is
$G(z)=z^{-1}\ln \bigl(1/(1-z)\bigr)$; therefore the convolution
family for $G(z)^x$ is $G_n(x)=xF_n(x+n)/(x+n)=x\sigma_n(x+n)$.

In this example the convolution family for
$e^{xzG(z)}=(1-z)^{-x}$ is 
$${x+n-1\choose n}={1\over n!}\left({n\brack 0}+{n\brack 1}\,x
+\,\cdots\, +{n\brack n}\,x^n\right)\,;$$
 therefore
$${n\brack n-k}={n!\over (n-k)!}\;G_k(n-k)={n!\over
(n-k)!}\,(n-k)\,\sigma_k(n)=n(n-1)\,\ldots\,(n-k)\,\sigma_k(n)\,.$$
We also have
$${n\brace n-k}={k-n\brack
-n}=(k-n)(k-1-n)\,\ldots\,(-n)\,\sigma_k(k-n)\,.$$
These formulas, which are polynomials in $n$ of degree~$2k$ for every
fixed~$k$, explain why the $\sigma$~functions are called Stirling
polynomials.
Notice that $\sigma_n(1)=(-1)^nB_n/n!$; it can also be shown that
$\sigma_n(0)=-B_n/(n\cdot n!)$.

The process of going from $F_n(x)$ to $xF_n(x+n)/(x+n)$ can be
iterated: Another replacement gives $xF_n(x+2n)/(x+2n)$, and after
$t$~iterations we discover that the polynomials $xF_n(x+tn)/(x+tn)$
also form a convolution family. This holds for all positive
integers~$t$, and the convolution condition is expressible as a set of
polynomial relations in~$t$; {\sl therefore\/
$xF_n(x+tn)/(x+tn)$ is a convolution family for all complex
numbers\/~$t$}. If $F_n(x)=[z^n]\,F(z)^x$, then
$xF_n(x+tn)/(x+tn)=[z^n]\,{\cal F}_t(z)^x$, where ${\cal F}_t(z)$ is
defined implicitly by the equation
$${\cal F}_t(z)=F\bigl(z{\cal F}_t(z)^t\bigr)\,.$$
In particular, we could have deduced the convolution
properties of the coefficients of ${\cal B}_t(z)^x$ in this way.

Let us restate what we have just proved, combining it with the
``derived convolution formula'' obtained earlier:

\proclaim Theorem. Let\/ $F_n(x)$ be any family of polynomials
in\/~$x$ such that\/ $F_n(x)$ has degree\/ $\leq n$. If
$$F_n(2x)=\sum_{k=0}^nF_k(x)\,F_{n-k}(x)$$
holds for all $n$ and $x$, then the following identities hold for all
$n$, $x$, $y$, and~$t$:
$$\eqalign{{(x+y)\,F_n(x+y+tn)\over x+y+tn}
&=\sum_{k=0}^n\; {x\,F_k(x+tk)\over
x+tk}\;{y\,F_{n-k}\bigl(y+t(n-k)\bigr)\over y+t(n-k)}\,;\cr
\noalign{\smallskip}
{n\,F_n(x+y+tn)\over x+y+tn}
&=\sum_{k=1}^n\;{k\,F_k(x+tk)\over
x+tk}\;{y\,F_{n-k}\bigl(y+t(n-k)\bigr)\over y+t(n-k)}\,.\cr}$$

\bn
{\bf Additional constructions.}\enspace
We have considered several ways to create new convolution families
from given ones, by multiplication or exponentiation of the associated
convolution matrices, or by replacing $F_n(x)$ by
$x\,F_n(x+tn)/(x+tn)$. It is also clear that the polynomials
$\alpha^nF_n(\beta x)$ form a convolution family whenever the
polynomials $F_n(x)$ do.

One further operation deserves to be mentioned: If $F_n(x)$ and
$G_n(x)$ are convolution families, then so is the family $H_n(x)$
defined by
$$H_n(x)=\sum_{k=0}^n\,F_k(x)\,G_{n-k}(x)\,.$$
This is obvious, since $H_n(x)=[z^n]\,F(z)^xG(z)^x$.
The corresponding operation on matrices $F=(f_{nk})$, $G=(g_{nk})$,
$H=(h_{nk})$ is 
$$h_{nk}=\sum_{i,j}\,{n\choose j}f_{ji}\,g_{(n-j)(k-i)}\,.$$
If we denote this binary operation by $H=F\circ G$, it is interesting
to observe that the associative law holds: $(E\circ F)\circ
G=E\circ(F\circ G)$ is true for all matrices $E$, $F$, $G$, not just
for convolution matrices. A~convolution matrix is characterized by
the special property $F\circ F=F\,{\rm diag}(2,4,8,\ldots\,)$.

The construction just mentioned is merely a special case of the
one-parameter family
$$H_n^{(t)}(x)=\sum_{k=0}^n\,F_k(x)\,G_{n-k}(x+tk)\,.$$
Again, $\{H_n^{(t)}(x)\}$ turns out to be a convolution family, for
arbitrary~$t$: We have
$$\eqalign{\sum_{k=0}^n H_n^{(t)}(x)z^n=\sum_{n\ge
k\ge0}\!\!F_k(x)\,G_{n-k}(x+tk)z^n
&=\sum_{n,k\ge0}\!F_k(x)\,G_n(x+tk)z^{n+k}\cr
&=\sum_{k\ge0}F_k(x)z^kG(z)^{x+tk}=G(z)^xF\bigl(zG(z)^t\bigr)^x,\cr}$$
so $H_n(x)=[z^n]\,\bigl(G(z)F\bigl(zG(z)^t\bigr)\bigr)^x$.
\bn
{\bf Applications.}\enspace
What's the use of all this? Well, we have shown that many interesting
convolution families exist, and that we can deduce nonobvious facts
with comparatively little effort once we know that we're dealing with
a convolution family.

One moral to be drawn is therefore the following. Whenever you encounter a
triangular pattern of numbers that you haven't seen before, check to
see if the first three rows have the form
$$\vcenter{\halign{\hfil$#$\hfil\quad&\hfil$#$\hfil\quad&\hfil$#$\hfil\cr
a\cr
b&a^2\cr
c&3ab&a^3\cr}}$$
for some $a,b,c$. (You may have to multiply or divide the $n$th row by
$n!$ first, and/or reflect its entries left to right.) If so, and if
the problem you are investigating is mathematically ``clean,'' chances
are good that the fourth row will look like
$$d\quad 4ac+3b^2\quad 6a^2b\quad a^4\,.$$
And if so, chances are excellent that you are dealing with a
convolution family. And if so, you may well be able to solve your
problem.

In fact, exactly that scenario has helped the author on several
occasions.

\bn
{\bf Asymptotics.}\enspace
Once you have identified a convolution family $F_n(x)$, you may well
want to know the approximate value of $F_n(x)$ when $n$ and~$x$ are
large. The remainder of this
paper discusses a remarkable general power series expansion, discovered
with the help of {\sl Mathematica}, which accounts for the behavior of
$F_n(x)$ when $n/x$ stays bounded and reasonably small as
$x\ra\infty$, although $n$ may also vary as a function of~$x$. 
We will assume that $F_n(x)$ is the coefficient of~$z^n$
in $F(z)^x$, where $F(0)=F'(0)=1$. 

Our starting point is the classical ``saddle point method,'' which
shows that in many cases the coefficient of~$z^n$ in a power series
$P(z)$ can be approximated by considering the value of~$P$ at a point
where the derivative of $P(z)/z^n$ is zero. (See [Good 1957].) In our
case we have $P(z)=e^{xf(z)}$, where $f(z)=\ln
F(z)=z+f_2z^2\!/2!+\cdots\,$; and the derivative is zero when
$x\,f'(z)=n/z$. Let this saddle point occur at $z=s$; thus, we have
$$s\,f'(s)=n/x\,.$$
Near $s$ we have $f(z)=f(s)+(z-s)f'(s)+O\bigl((z-s)^2\bigr)$; so we
will base our approximation on the assumption that the
$O\bigl((z-s)^2\bigr)$ contribution is zero. The approximation to $F_n(x)$
will be $\widetilde{F}_n(x)$, where
$$\eqalign{%
\widetilde{F}_n(x)&=[z^n]\,\exp\,\bigl(x\,f(s)+x\,(z-s)\,f'(s)\bigr)\cr
\noalign{\medskip}
&={e^{x(f(s)-sf'(s))}\over n!}\;x^nf'(s)^n={F(s)^x\over
n!}\;\left({n\over es}\right)^n\,.\cr}$$

First let's look at some examples; later we will show that the ratio
$F_n(x)/\widetilde{F}_n(x)$ is well behaved as a formal power series.
Throughout this discussion we will let
$$y=n/x\,;$$
our goal, remember, is to find approximations that are valid when $y$
is not too large, as $x$ and possibly~$n$ go to~$\infty$.

The simplest example is, of course, $F(z)=e^z$ and $f(z)=z$; but we
needn't sneeze at it because it will give us some useful calibration.
In this case $f''(z)=0$, so our approximation will be exact. We have
$s=y$, hence
$$\widetilde{F}_n(x)={e^{xy}\over n!}\,\left({n\over
ey}\right)^n={e^n\over n!}\,\left({x\over e}\right)^n={x^n\over
n!}=F_n(x)\,.$$

Next let's consider the case $F(z)=T(z)/z$, $f(z)=T(z)$, when we know
that $F_n(x)=x(x+n)^{n-1}/n!$. In this case
$z\,T'(z)=T(z)/\bigl(1-T(z)\bigr)$, so we have
$T(s)/\bigl(1-T(s)\bigr)=y$ or
$$T(s)={y\over 1+y}\,,\qquad s={y\over 1+y}\,e^{-y/(1+y)}$$
because $T(z)=ze^{T(z)}$. Therefore
$$\widetilde{F}_n(x)={e^{xy/(1+y)}\over n!}\,\left({n(1+y)\over
ey\,e^{-y/(1+y)}}\right)^n={(x+n)^n\over n!}\,;$$
the ratio $F_n(x)/\widetilde{F}_n(x)=x/(x+n)=1/(1+y)$ is indeed near~1
when $y$ is small.

If $F(z)=1+z$ we find, similarly, $s=y/(1-y)$ and
$$n!\,\widetilde{F}_n(x)=\left({1\over
1-y}\right)^x\,\left({n(1-y)\over ey}\right)^n={x^xe^{-n}\over
(x-n)^{x-n}}\,;$$ 
by Stirling's approximation we also have
$$n!\,F_n(x)={x!\over (x-n)!}={x^xe^{-n}\over
(x-n)^{x-n}}\,(1-y)^{-1/2}\bigl(1+O(x^{-1})\bigr)\,.$$
Again the ratio $F_n(x)/\widetilde{F}_n(x)$ is near~1.
In general if $F(z)={\cal B}_t(z)$ the saddle point~$s$ turns out to be
$y\bigl(1+(t-1)y\bigr)^{t-1}/(1+ty)^t$, and
$$n!\,\widetilde{F}_n(x)={(x+tn)^{x+tn}e^{-n}\over
\bigl(x+(t-1)n\bigr)^{x+(t-1)n}}\,;$$
a similar analysis shows that this approximation is quite good, for
any fixed~$t$.

We know that
$$F_n(x)={x^n\over n!}\,\left(1+{f_{n(n-1)}\over x}+{f_{n(n-2)}\over
x^2}+\cdots\,\right)$$
and that $f_{n(n-k)}$ is always a polynomial in $n$ of degree $\leq 2k$.
Therefore if $n^2\!/x\ra 0$ as $x\ra\infty$, we can simply use the
approximation $F_n(x)=(x^n\!/n!)\bigl(1+O(n^2\!/x)\bigr)$. But there are
many applications where we need a good estimate of $F_n(x)$ when
$n^2\!/x\ra\infty$ while $n/x\ra 0$; for example, $x$~might be
$n\log n$. In such cases $\widetilde{F}_n(x)$ is close to $F_n(x)$ but
$x^n\!/n!$ is not. 

We can express $s/y$ as a power series in $y$ by inverting the power
series expression $sf'(s)=y$:
$$s/y=1-f_2y+(4f_2^2-f_3)y^2\!/2+(15f_2f_3-30f_2^3-f_4)y^3\!/6+\cdots\;.$$
From this we can get a formal series for $\widetilde{F}_n(x)$,
$$\eqalign{\widetilde{F}_n(x)
&={x^n\over
n!}\;{\exp\bigl(n(s/y)(1+f_2s/2!+f_3s^2\!/3!+\cdots\,)-n\bigr)\over 
(s/y)^n}\cr
\noalign{\medskip}
&={x^n\over n!}\;\left(1+{nf_2\over
2}\;y+{3n^2f_2^2-12nf_2^2+4nf_3\over
24}\;y^2+O(n^3y^3)\right)\,.\cr}$$
We can also use the formula
$$f_{n(n-k)}=\sum\;{n^{\underline{k+k_2+k_3+\cdots}}\over
2!^{k_2}\,k_2!\,3!^{k_3}\,k_3!\,\ldots}\;f_2^{k_2}f_3^{k_3}\,\ldots\;,$$
where the sum is over all nonnegative $k_2,k_3,\ldots$ with
$k_2+2k_3+\cdots =k$, to write
$$\eqalign{F_n(x)={x^n\over n!}\;&\left(1+{nf_2-f_2+O(x^{-1})\over
2}\;y\right.\cr
\noalign{\smallskip}
&\qquad \null+\left.{3n^2f_2^2-18nf_2^2+4nf_3+33f_2^2-12f_3+O(x^{-1})\over
24}\;y^2+O(n^3y^3)\right)\,.\cr}$$
These series are not useful asymptotically unless $ny=n^2\!/x$ is
small. But the approximation 
$\widetilde{F}_n(x)$ itself is excellent, because amazing
cancellations occur when we compute the ratio:
$${F_n(x)\over\widetilde{F}_n(x)}=1-{f_2\over 2}\;y+{11f_2^2-4f_3\over
8}\;y^2+O(y^3)+O(x^{-1})\,.$$

\proclaim Theorem. When
$F(z)=\exp(z+f_2z^2\!/2!+f_3z^3\!/3!+\cdots\,)$ and the functions
$F_n(x)$ and
$\widetilde{F}_n(x)$ are defined as above, the ratio
$F_n(x)/\widetilde{F}_n(x)$ can be written as a formal power series
$\sum_{i,j\geq 0}c_{ij}y^ix^{-j}$, where $y=n/x$ and the
coefficients~$c_{ij}$ are polynomials in $f_2,f_3,\ldots\;$.

The derivation just given shows that we can write
$F_n(x)/\widetilde{F}_n(x)$ as a formal power series of the form
$\sum_{i,j\geq 0}a_{ij}n^ix^{-j}$, where $a_{ij}=0$ when $i>2j$; the
surprising thing is that we also have $a_{ij}=0$ whenever $i>j$.
Therefore we can let $c_{ij}=a_{i(i+j)}$.

To prove the theorem, we let $R(z)=1+R_1z+R_2z^2+\cdots$ stand for the
terms neglected in our approximation:
$$F(z)^x=e^{xf(s)-xsf'(s)}\left(1+{n\over s}\;{z\over 1!}+{n^2\over
s^2}\;{z^2\over 2!}+{n^3\over s^3}\;{z^3\over
3!}+\cdots\,\right)\,R(z)\,.$$
The coefficient of $z^n$ is
$$F_n(x)=\widetilde{F}_n(x)\left(1+R_1s+{n-1\over
n}\;R_2s^2+{(n-1)(n-2)\over n^2}\;R_3s^3+\cdots\,\right)\,;$$
so the ratio $F_n(x)/\widetilde{F}_n(x)$ is equal to
$$\sum_{k\geq 0}\,{n^{\underline{k}}\over n^k}\,R_ks^k=\sum_{j,k\geq
0}(-n)^{-j}\,{k\brack k-j}\,R_ks^k=\sum_j(-n)^{-j}P_j\,,$$
where $P_j=\sum_k{k\brack k-j}R_ks^k$ is a certain power series in~$s$
and~$x$. The coefficients~$R_k$ are themselves power series in~$s$
and~$x$, because we have
$$R(z)=\exp\left(x(z-s)^2\;{f''(s)\over 2!}+x(z-s)^3\;{f'''(s)\over
3!}+\cdots\,\right)\,.$$

We know from the discussion above that
$${k\brack k-j}=k(k-1)\,\ldots\,(k-j)\,\sigma_j(k)$$
is a polynomial in $k$. Therefore we can write
$$P_j=\left.{\vartheta\brack\vartheta-j}\,R(z)\,\right\vert_{z=s}\,,$$
where $\vartheta$ is the operator that takes $z^k\mapsto k\,z^k$ for
all~$k$; i.e., $\vartheta G(z)=z\,G'(z)$ for all power series $G(z)$.
The theorem will be proved if we can show that $P_j/n^j$ is a formal
power series in~$y$ and~$x^{-1}$, and if the sum of these formal power
series over all~$j$ is also such a series.

Consider, for example, the simplest case $P_0=R(s)$; obviously
$P_0=1$. The next simplest case is
$P_1=\left.{\vartheta\brack\vartheta-1}\,R(z)\,\right\vert_{z=s}=
\left.{1\over 2}\vartheta(\vartheta-1)R(z)\,\right\vert_{z=s}$. It is
easy to see that
$$\vartheta^{\underline{j}}=z^jD^j\,,$$
where $D$ is the differentiation operator $D\,G(z)=G'(z)$, because
$z^jD^j$ takes $z^k$ into $k^{\underline{j}}z^k$. Therefore
$$P_1={\textstyle{1\over 2}} s^2R''(s)={\textstyle{1\over
2}}\,xs^2f''(s)\,.$$ 
It follows that $P_1/n={1\over 2}(s/y)sf''(s)$ is a power series
in~$y$; it begins ${1\over 2}f_2y+{1\over 2}(f_3-f_2^2)y^2+\cdots\;$.

Now let's consider $P_j$ in general. We will use the fact that the
Stirling numbers ${k\brack k-j}$ can be represented in the form
$${k\brack k-j}=p_{j1}{k\choose j+1}+p_{j2}{k\choose j+2}+\cdots
+p_{jj}{k\choose 2j}\,,$$
where the coefficients $p_{ji}$ are the positive integers in the
following triangular array:
$$\vcenter{\halign{\hfil$#$\quad&\hfil$#$\quad&\hfil$#$\quad
&\hfil$#$\quad&\hfil$#$\quad&\hfil$#$\cr
1\cr
2&3\cr
6&20&15\cr
24&130&210&105\cr
120&924&2380&2520&945&.\cr}}$$
$\bigl($This array is clearly not a convolution matrix; but the theory
developed above implies that the numbers $j!\,p_{ji}/(i+j)!\,$, namely
$$\vcenter{\halign{\hfil$#$\hfil\quad&\hfil$#$\hfil\quad&\hfil$#$\hfil\quad
&\hfil$#$\hfil\quad&\hfil$#$\hfil\quad&\hfil$#$\hfil\cr
1/2\cr
2/3&1/4\cr
3/2&1&1/8&&&,\cr
24/5&13/3&1&1/16\cr
20&22&85/12&5/6&1/32\cr}}$$
do form the convolution matrix for the powers of 
$\exp(z/2+z^2\!/3+z^3\!/4+\cdots\,)$.
The expression ${k\brack k-j}=\sum_{i=1}^jp_{ji}{k\choose j+i}$ was
independently discovered by [Appell 1880], [Jordan 1933], and [Ward
1934]. The number of permutations of $i+j$ elements having no fixed
points and exactly $i$~cycles is~$p_{ji}$, an ``associated Stirling
number of the first kind'' [Riordan 1958, section 4.4] 
[Comtet 1970, exercise 6.7].$\bigr)$
It follows that
$$P_j=p_{j1}s^{j+1}\;{R^{(j+1)}(s)\over
(j+1)!}+p_{j2}s^{j+2}\;{R^{(j+2)}(s)\over
(j+2)!}+\cdots+p_{jj}s^{2j}\;{R^{(2j)}(s)\over (2j)!}\,.$$

Now $R(z)$ is a sum of terms having the form
$$a_{il}x^i(z-s)^l\,,$$
where $l\geq 2i$ and where $a_{il}$ is a power series in~$s$. Such a
term contributes $a_{il}x^is^lp_{j(l-j)}$ to~$P_j$; so it contributes
$a_{il}(s/y)^js^{l-j}x^{i-j}p_{j(l-j)}$ to $P_j/n^j$. This contribution
is nonzero only if $j<l\leq 2j$. Since $l\geq 2i$, we have $i\leq j$;
so $P_j/n^j$ is a power series in~$y$ and~$x^{-1}$.

For a fixed value of $j-i$, the smallest power of $y$ that can occur in
$P_j/n^j$ is $y^{2i-j}=y^{j-2(j-i)}$. Therefore only a finite number
of terms of $\sum_jP_j/(-n)^j$ contribute to any given power of~$y$
and~$x^{-1}$. This completes the proof.

A careful analysis of the proof, and a bit of {\sl Mathematica\/}
hacking, yields the more precise result
$${F_n(x)\over \widetilde{F}_n(x)}={1\over
(1+s^2y^{-1}d_2)^{1/2}}+{(s/y)^3A\over
x(1+s^2y^{-1}d_2)^{7/2}}+O(x^{-2})\,,$$
where $A={1\over 12}s^3y^{-1}d_2^3-{3\over 4}sd_2^2-{1\over
2}s^2d_2d_3-{5\over 24}s^3d_3^2+{1\over 3}yd_3+{1\over 8}s^3d_2d_4+
{1\over 8}syd_4$ and $d_k=f^{(k)}(s)$.

\bn
{\bf Acknowledgment.}\enspace
I wish to thank Ira Gessel and Svante Janson for stimulating my
interest in this subject and for their helpful comments on the first
draft. Ira Gessel and Richard Brent also introduced me to several
relevant references.

\bn
{\bf References}
\def\bibspace{\unskip\quad}

\bib
Appell, P. \bibspace 1880.\bibspace
``D\'eveloppement en s\'erie enti\`ere de $(1+ax)^{1/x}$.'' {\sl
Archiv der Mathematik und Physik\/ \bf 65}: 171--175.

\bib
Bell, E. T. \bibspace 1934. \bibspace ``Exponential numbers.''
{\sl American Mathematical Monthly\/ \bf 41}: 411--419.


\bib
Carlitz, L.\bibspace 1978.\bibspace
``Generalized Stirling and related numbers.'' {\sl Rivista di
Matematica della Universit\`a di Parma}, serie~4, {\bf 4}: 79--99.

\bib
Cayley, A. \bibspace 1859. \bibspace
 ``On the analytical forms called trees. Second
part.''
{\sl Philosophical Magazine\/ \bf 18}: 371--378.
Reprinted in Cayley's
{\sl Collected Mathematical Papers\/ \bf 4}: 112--115.

\bib
Comtet, Louis. \bibspace 1970. \bibspace
 {\sl Analyse Combinatoire}. Presses
Universitaires de France. (English translation, {\sl Advanced
Combinatorics}, D.~Reidel, Dordrecht, 1974.)

\bib
Garsia, Adriano M.\bibspace 1973.\bibspace
``An expos\'e of the Mullin-Rota theory of polynomials of binomial
type.'' {\sl Linear and Multilinear Algebra\/ \bf 1}: 47--65.

\bib
Gomes Teixeira, F. \bibspace 1896. \bibspace
``Sur le d\'eveloppement de $x^k$ en s\'erie ordonn\'ee suivant les
puissances du sinus de la variable.'' {\sl Nouvelles Annales de
Math\'ematiques}, s\'erie 3, {\bf 15}: 270--274.

\bib
Good, I. J. \bibspace 1957. \bibspace
 ``Saddle-point methods for the multinomial
distribution.'' {\sl Annals of Mathematical Statistics\/ \bf 28}:
861--881.

\bib
Gould, H. W., and Kauck\'y, I. \bibspace 1966. \bibspace
 ``Evaluation of a class of
binomial coefficient summations.'' {\sl Journal of Combinatorial
Theory\/ \bf 1}: 233--248.

\bib
Graham, Ronald L., Knuth, Donald E., and Patashnik, Oren. \bibspace
 1989. \bibspace
{\sl Concrete Mathematics}. Addison-Wesley, Reading, Massachusetts.

\bib
Harris, Bernard,  and Schoenfeld, Lowell. \bibspace 1967. \bibspace
 ``The number of
idempotent elements in symmetric semigroups.'' {\sl Journal of
Combinatorial Theory\/ \bf 3}: 122--135.

\bib
Jabotinsky, Eri. \bibspace  1947. \bibspace
 ``Sur la repr\'esentation
de la composition de fonctions par un produit de matrices. Applicaton
\`a l'it\'eration de~$e^x$ et de $e^x-1$.'' {\sl Comptes Rendus
Hebdomadaires des Sciences de L'Academie des Sciences}, {\bf 224}: 323--324.

\bib
Jabotinsky, Eri. \bibspace 1953. \bibspace
``Representation of functions by matrices.
Application to Faber polynomials.'' {\sl Proceedings of the American
Mathematical Society\/ \bf 4}: 546--553.

\bib
Jabotinsky, Eri. \bibspace 1963. \bibspace
 ``Analytic iteration.'' {\sl Transactions of
the American Mathematical Society\/ \bf 108}: 457--477.

\bib
Jordan, Charles.\bibspace 1933.\bibspace
``On Stirling's numbers.'' {\sl T{\^o}hoku Mathematical Journal\/ \bf
37}: 254--278.

\bib
Knuth, Donald E. \bibspace 1992. \bibspace
 ``Two notes on notation.'' {\sl American
Mathematical Monthly\/ \bf 99}: 403--422.

\bib
Knuth, Donald E. and Pittel, Boris. \bibspace 1989. \bibspace
 ``A~recurrence related to
trees.''
{\sl Proceedings of the American Mathematical Society\/ \bf 105}:
335--349.

\bib
Kuczma, Marek.\bibspace 1968.\bibspace
{\sl Functional Equations in a Single Variable}. Polish Scientific
Publishers, Warsaw.

\bib
Lah, I. \bibspace  1955. \bibspace
 ``Eine neue Art von Zahlen, ihre Eigenschaften und
Anwendung in der mathemat\-ischen Statistik,'' {\sl Mitteilungsblatt
f\"ur Mathematische Statistik\/ \bf 7}: 203--212.

\bib
Motzkin, T. S. \bibspace  1971. \bibspace
 ``Sorting numbers for cylinders and other
classification numbers.'' {\sl Proceedings of Symposia in Pure
Mathematics\/ \bf 19}: 167--176.

\bib
Mullin, Ronald, and Rota, Gian-Carlo. \bibspace 1970.  \bibspace
``On the foundations of combinatorial theory. III. Theory of binomial
enumeration.'' In {\sl Graph Theory and Its Applications}, edited by
Bernard Harris (Academic Press, 1970), 167--213.

\bib
Riordan, John.\bibspace 1958.\bibspace
{\sl An Introduction to Combinatorial Analysis}. John Wiley \& Sons,
New York.

\bib
Riordan, John. \bibspace  1968. \bibspace
 {\sl Combinatorial Identities}. John Wiley \&
Sons, New York.

\bib
Rogers, D. G.\bibspace 1978.\bibspace
``Pascal triangles, Catalan numbers and renewal arrays.'' {\sl
Discrete Mathematics\/ \bf 22}: 301--310.

\bib
Roman, Steven.\bibspace 1978.\bibspace
{\sl The Umbral Calculus}. Academic Press, Orlando.

\bib
Roman, Steven M., and Rota, Gian-Carlo. \bibspace 1978. \bibspace
``The umbral calculus.'' {\sl Advances in Mathematics\/ \bf 27}:
95--188.

\bib
Rota, Gian-Carlo, Kahaner, D., and Odlyzko, A. \bibspace 1973. \bibspace
``On the foundations of combinatorial theory. VIII. Finite operator
calculus.'' {\sl Journal of Mathematical Analysis and Applications\/ \bf
42}: 884--760. Reprinted in Rota, Gian-Carlo, {\sl Finite Operator
Calculus\/} (Academic Press, 1975), 7--82.

\bib
Steffensen, J. F.\bibspace 1941.\bibspace
``The poweroid, an extension of the mathematical notion of power.''
{\sl Acta Mathematica\/ \bf 73}: 333--366.

\bib
Ward, Morgan.\bibspace 1934.\bibspace
``The representation of Stirling's numbers and Stirling's polynomials
as sums of factorials.'' {\sl American Journal of Mathematics\/ \bf
56}: 87--95.

\bye